\documentclass[12pt]{amsart}
 
\usepackage{amssymb,amsfonts}
\usepackage[all,arc]{xy}
\usepackage{enumerate}
\usepackage{mathrsfs}
\usepackage{enumitem}

\newcommand{\norm}[1]{\left\lVert#1\right\rVert}
\usepackage[
margin=2.6cm,
includefoot,
footskip=20pt,
]{geometry}
\usepackage{hyperref}
\newtheorem{thm}{Theorem}[section]

\newtheorem{prop}[thm]{Proposition}

\theoremstyle{definition}

\theoremstyle{remark}
\newtheorem{rem}[thm]{Remark}

\makeatletter
\@addtoreset{equation}{section}
\makeatother
\numberwithin{equation}{section}
\begin{document}
	
	\setlength{\baselineskip}{18truept}
	
	\begin{center}
		
		{\Large\bf An unstable three dimensional KAM torus for the quintic NLS}\\
		\vspace*{1cm}
		
				\bfseries
					 {Nguyen Thuy Trung}\\
		\vspace*{1cm}
		{\bf NANTES - 2018 }\\
	\end{center}
	
	\thispagestyle{empty}
	
	\title{}

	\maketitle

	\tableofcontents
		Abstract: 
	  We consider the quintic nonlinear Schr\"{o}dinger on the circle. By applying a Birkhoff procedure and a KAM theorem, we exihibit a three dimension invariant torus that is linearly unstable. In comparison, we also prove that two dimensional tori are always linearly stable.     
	\section{Introduction}
We consider the non linear Schr\"{o}dinger equation on the torus 
\begin{equation}
i \partial_t u + \partial_{ x x} u = \vert u^4 \vert u, \quad ( t, x) \in \mathbb{R} \times \mathbb{T}. \label{intro1}
\end{equation} 
This is an infinite dimensional dynamic system on the phase space $ (u, \bar{u}) \in  L^2(\mathbb{T})$
endowed with the symplectic form $ -idu \wedge d \bar{u}.$ 
The flow $u(t)$ preserves the Hamiltonian
$$ h = \int_{\mathbb{T}} \vert u_x \vert^2 + \frac{1}{3}\vert u \vert^6 dx,$$ 
and also, the mass and the momentum 
         $$ \mathbb{L} = \int_{\mathbb{T}} |u|^2 dx, \quad \, \mathbb{M} = \int_{\mathbb{T}} Im ( u \cdot \nabla \bar{u})dx.$$
Let us expand $ u$ and $ \bar{u}$ in Fourier basis: 
$$ u(t,x) = \sum_{j \in \mathbb{Z} } a_j (t) e^{ i j x}, \quad \bar{u}(t,x) = \sum_{j \in \mathbb{Z} } b_j(t) e^{ -i j x}.$$
In this variables, the symplectic structure becomes $ -i \sum_{j \in \mathbb{Z} } d a_j \wedge d b_j.$ The Hamiltonian h of the system reads
 $$ h =  \sum_{j \in \mathbb{Z}} j^2 a_j b_j + \frac{1}{3}\sum_{j,\ell \in \mathbb{Z}^3; \mathcal{M}(j,l) = 0 } a_{j_1}a_{j_2}a_{j_3}b_{\ell_1}b_{\ell_2}b_{\ell_3}=N+P,$$
 and the mass and the momentum
           $$ \mathbb{L} = \sum_{j \in \mathbb{Z}} a_j b_j, \quad \, \mathbb{M} = \sum_{j \in \mathbb{Z}} j a_j b_j, $$
 here $ \mathcal{M}(j,l) = j_1 + j_2 +j_3 - \ell_1 - \ell_2 -\ell_3$ denotes the momentum of the monomial $a_{j_1}a_{j_2}a_{j_3}b_{\ell_1}b_{\ell_2}b_{\ell_3}.$ We can rewrite equation \eqref{intro1} into a system of infinite equations 
                       $$ \begin{cases}
                       i \dot{a_j} &= j^2 a_j + \frac{\partial P}{ \partial b_j} \quad j \in \mathbb{Z}, \\
                       -i \dot{b_j} &= j^2 b_j + \frac{\partial P}{ \partial a_j} \quad j \in \mathbb{Z}.
                       \end{cases}	$$

 In this article, we are interesting in the dynamic behavior near to $0$ of solution of \eqref{intro1} in two specific forms:
             \begin{equation}
             u(t,x) =  a_p (t) e^{ipx} e^{-ip^2 t} +  a_q (t) e^{iqx} e^{-iq^2 t}  + \mathcal{O} (\varepsilon), \label{upq} 
             \end{equation} 
        and
             \begin{equation}
             u(t,x) =  a_p (t) e^{ipx} e^{-ip^2 t}  +  a_q (t) e^{iqx} e^{-iq^2 t} + a_m (t) e^{imx} e^{-im^2 t}  + \mathcal{O} (\varepsilon), \label{upqm} 
             \end{equation} 
 or more precisely the persistence of two and three dimensional linear
 invariant tori:
              \begin{align}
             	 \mathbf{T}_c^{2}(p,q) &= \{ | a_p|^2 =c_1,\; | a_q|^2 = c_2\}, \label{2dimtorus}\\
             	  \mathbf{T}_c^{3}(p,q,m) &= \{ | a_p|^2 =c_1,\; | a_q|^2 = c_2,\;| a_m|^2 = c_3  \}\label{3dimtorus},
             \end{align}	
      with $0 < c_1,c_2,c_3 \ll 1.$\\
 The first result of this paper is stated for two dimensional tori.
 \begin{thm} \label{thm2}
 	Fix $ p, \,q \in \mathbb{Z},$ and $ s > \frac{1}{2}.$ There exists $\nu_0 >0,$ and for $ 0 < \nu < \nu_0,$ there exists $\mathcal{ D}_{\nu} \subset [ 1,2]^2$ asymptotically of full measure \textup{(}i.e. $meas ([ 1,2]^2 \setminus \mathcal{ D}_{\nu}) \to 0$ when $ \nu \to 0$\textup{)} such that for $ \rho \in \mathcal{ D}_{\nu},$ equation \eqref{intro1} admits a solution of the form 
 	$$  u ( x) = \sum_{ j \in \mathbb{Z} } a_j ( t \omega) e^{ijx}$$
 	where $ \{a_j\}_j $ is analytic function from $ \mathbb{T}^2 $ to $ \ell_{s}^2$ satisfying uniformly in $ \theta \in \mathbb{T}^2$
 	$$ | a_p - \sqrt{\nu \rho_1}|^2 + | a_q - \sqrt{\nu \rho_2}|^2 + \sum_{ j \neq p,q} (1 + j^2)^s | a_j|^2 = \mathcal{ O} ( \nu^2).$$
 	Here $\omega $ is a nonresonant vector in $\mathbb{R}^2$ that satisfies 
 	$$ \omega = ( p^2, q^2) + \mathcal{ O}(\nu^2).$$ 
 	Furthermore, this solution is linearly stable. 
 \end{thm} 
  For three dimensional tori, it is too complicated\footnote{the difficulty is to verify KAM hypotheses} to consider the general case.
 In order to apply KAM theorem \ref{mainthm}, we avoid the case where there is $\ell \in \mathbb{Z} $ solving equation \footnote{ in this case, the linear part $a_{j_1}^2a_{j_2}b_{j_3}^2b_{\ell} + b_{j_1}^2 b_{j_2}a_{j_3}^2 a_{\ell}$ of the mode $\ell$ would create the instability, and the energy would soon transfer mainly between four modes $p,q,m,\ell$, which was studied carefully in \cite{quinticNLSeq}.   }
 \begin{equation}
 \label{introell} \begin{cases}
 2j_1 + j_2 &= 2j_3 + \ell  \\
 2j_1^2 +j_2^2 &=2j_3^2+ \ell^2.
 \end{cases}
  \end{equation}
 In this paper, we will give here an example of $(p,q,m)$ and $\rho$ such that for $ \nu$ small enough the torus $ \mathbf{T}^{3}_{\nu \rho}(p,q,m) = \{ | a_p|^2 = \nu \rho_1,\; | a_q|^2 = \nu \rho_2, \; | a_m|^2 = \nu \rho_3\}$ is linearly unstable. For $\epsilon = 10^{-2}$, denote 
  $$\mathcal{D}=\mathcal{D}_2 = [2-\epsilon, 2+\epsilon] \times[1-\epsilon, 1+\epsilon] \times[9-\epsilon, 9+\epsilon] . $$

    \begin{thm} \label{thm3}
    	Fix $ p=-3, \,q=10,\,m=-6,$ and $ s > \frac{1}{2}$. There exists $\nu_0 >0,$ and for $ 0 < \nu < \nu_0,$ there exists $\mathcal{ D}_{\nu} \subset \mathcal{ D}$ asymptotically of full measure \textup{(}i.e. $meas (\mathcal{ D} \setminus \mathcal{ D}_{\nu}) \to 0$ when $ \nu \to 0$\textup{)} such that for $ \rho \in \mathcal{ D}_{\nu},$ equation \eqref{intro1} admits a solution of the form 
    	\begin{align} 
    	u ( x) = \sum_{ j \in \mathbb{Z} } a_j ( t \omega) e^{ijx}
    	\end{align}
    	where $ \{a_j\}_j $ is analytic function from $ \mathbb{T}^3 $ to $ \ell_{s}^2$ satisfying uniformly in $ \theta \in \mathbb{T}^3$
    	\begin{align} 
    	| a_p - \sqrt{\nu \rho_1}|^2 + | a_q - \sqrt{\nu \rho_2}|^2 + | a_m - \sqrt{\nu \rho_3}|^2 + \sum_{ j \neq p,q,m} (1 + j^2)^s | a_j|^2 = \mathcal{ O} ( \nu^2).
    	\end{align} 
    	Here $\omega $ is a non resonant vector in $\mathbb{R}^3$ that satisfies 
    	$$ \omega = ( 3^2, 10^2, 6^2) + \mathcal{ O}(\nu^2).$$ 
    	Furthermore, this solution is linearly unstable.
    \end{thm} 
  In order to prove theorems \ref{thm2}, \ref{thm3}, we follow a general stratery developed in \cite{maintheorem} for a  system of coupled nonlinear Schr\"{o}dinger equations on the torus. Firstly, we apply a Birkhoff normal form procedure (Proposition \eqref{Birkhoff}) to kill the non resonances of $P$. Then we use sympletic changes of variables to diagonalize the effective part into the form of $h_0$. The hyperbolic directions of torus $\mathbf{T}^3_{\nu \rho}(-3,10,-6)$ are revealed in this step. Readers are suggested to take a look at the original statement of KAM theorem in \cite{maintheorem} for further understanding. 
  
The study of finite dimensional tori in an infinite dimensional phase space was pioneered by J. Bourgain \cite{Bourgain} in 1988. However, the existence of unstable KAM tori in one dimensional context was first proved by B. Gr\'{e}bert and V. Rocha \cite{maintheorem} in 2017, where they studied the system of coupled nonlinear Schr\"{o}dinger equations on the torus. For the equation \eqref{intro1}, in case of $u(0,x)$ supported maninly in four modes $(p,q,m,s)$, which satisfy such a relation in \eqref{introell}, the study of solutinon was studied carefully in \cite{quinticNLSeq} and \cite{H-P}. In particular, in \cite{H-P} they proved the recurrent exchange of energy between those modes.

\textit{Acknowledgement}: I wish to thank Professor Bernoit Gr\'{e}bert for motivating me to publish this paper with numerous suggestions and discussions. I also wish to thank Le Quoc Tuan and Lan Anh for computations in the appendix A. 

	\section{KAM theorem}

   In order to proof theorems \ref{thm2} and \ref{thm3}, we recall a KAM theorem stated in \cite{maintheorem}. \\
   We consider a Hamiltonian $h = h_0 + f,$ where $ h_0$ is a quadratic Hamiltonian in normal form 
 \begin{align}
 	h_0 = \Omega(\rho) \cdot r + \sum_{a \in \mathcal{Z} } \Lambda_a (\rho) \vert \zeta_a \vert^2.
 \end{align}	
Here 
 \begin{itemize}
 	\item $\rho$ is a parameter in $ \mathcal{D},$ which is a compact in the space $\mathbb{R}^n;$
 	\item $r \in \mathbb{R}^n $ are the actions corresponding to the internal modes $ ( r, \theta) \in \left( \mathbb{R}^n \times \mathbb{T}^n ,    dr  \wedge d\theta\right);$
 	\item $\mathcal{L}$ and $ \mathcal{F}$ are respectively infinite and finite sets, $ \mathcal{Z}$ is the disjoint uninon $ \mathcal{L} \cup   \mathcal{F};$
 	\item $ \zeta = \left( \zeta_a \right)_{ a\in \mathcal{Z}} \in \mathbb{C}^{\mathcal{Z}}$ are the external modes endowed with the standard complex symplectic structure $ -i d\zeta \wedge d\eta.$ The external modes decomposes in a infinite part $ \zeta_{ \mathcal{L}} = \left( \zeta_a\right)_{ a \in \mathcal{L} },$ corresponding to elliptic directions, which means $ \Lambda_{a} \in \mathbb{R}$ for $ a \in \mathcal{ L}$, and a finite part $\zeta_{ \mathcal{F}} = \left( \zeta_a\right)_{ a \in \mathcal{F} }, $ corresponding to hyperbolic directions, which means $ \Im \Lambda_a \neq 0 $ for $ a \in \mathcal{F};$
 	\item $ \mathcal{ L}$ has a clustering structure $ \mathcal{ L} = \cup_{ j \in \mathbb{N}} \mathcal{ L}_j,$ where $ \mathcal{ L}_j$ are finite sets of cardinality $ d_j \le d < \infty.$ If $ a \in \mathcal{ L}_j,$ we denote $ [a] = \mathcal{ L}_j$ and $ w_a = j,$ for $ a \in \mathcal{F}$ we set $w_a = 1;$ 
 	\item  the mappings 
 	      \begin{align}
 	      	&\Omega : \, \mathcal{ D} \to \mathbb{R}^n,\\
 	      	&\Lambda_a : \, \mathcal{ D} \to \mathbb{C}, \quad a \in \mathcal{ Z}, 
 	      \end{align}
       are smooth;
     \item$f = f( r, \theta, \zeta ; \rho)$ is a perturbation, small compare to the integrable part $h_0.$  
 \end{itemize}

\textbf{Linear space} Let $ s \ge 0,$ we consider the complex weighted $ \ell_2-$ space 
        $$ Z_s = \{ \zeta = \left( \zeta_a \in \mathbb{C}, \, a \in \mathcal{ Z} \right) | \norm{ \zeta}_s < \infty \},$$
     where $$ \norm{ \zeta}_s = \sum_{a \in \mathcal{ Z}} \vert \zeta_a \vert^2 w_a^{2s}.$$
   Similarly we difine 
        $$ Y_s = \{ \zeta_{\mathcal{L}} = \left( \zeta_a \in \mathbb{C}, \, a \in \mathcal{ L} \right) | \norm{ \zeta_{ \mathcal{L}}}_s  < \infty \},$$
        with the same norm. 
    We endow $Z_s \times Z_s$ and $Y_s \times Y_s$ with the symplectic structure $ -i d \zeta \wedge d \eta,$ with $ \eta = \bar{ \zeta}.$

 \textbf{A class of Hamiltonian functions.}  Denote $ \omega = (\zeta  , \eta).$  On the space 
          $$ \mathbb{C}^n \times \mathbb{C}^n \times \left( Z_s \times Z_s\right)$$
     we define the norm 
         $$ \norm{ \left( r, \theta, \omega \right) }_s = \max \left( \vert r \vert , \vert \theta \vert, \norm{ \zeta}_s \right).$$
     For $ \sigma >0$ we denote 
         $$ \mathbb{T}^n_{\sigma} = \{ \theta \in \mathbb{C}^n: \vert \Im \theta \vert < \sigma \} / 2 \pi \mathbb{Z}^n.$$
     For $ \sigma, \mu \in \left( 0,1 \right]$ and $ s \ge 0$ we set 
         $$ \mathcal{O}^s ( \sigma, \mu) = \{ r \in \mathbb{C}^n : \vert r \vert < \mu^2 \} \times \mathbb{T}^n_s \times \{ \omega  \in Z_s \times Z_s : \norm{ \zeta}_s < \mu \}.$$
     We will denote points in $ \mathcal{O}^s ( \sigma, \mu)$ as $ x = ( r, \theta, \omega).$ Let $ f : \mathcal{O}^0( \sigma, \mu) \times \mathcal{ D} \to \mathbb{C}$ be a $C^1 \textendash$function\footnote{$C^1$ regularity with respect to $\rho$ in the Whitney sense}, real holomorphic in the first variable $x$, such that for all $ \rho \in \mathcal{ D}, \, x \in \mathcal{O}^s(\sigma, \mu):$
      $$ \nabla_{\omega} f( x, \rho) \in Z_s \times Z_s$$
      and 
      $$ \nabla^2_{\omega_{\mathcal{ L} } \omega_{\mathcal{ L}} } f( x, \rho) \in \mathcal{ L} ( Y_s, Y_s)$$
      are real holomorphic functions. We denote by $ \mathcal{T}^s ( \sigma, \mu, \mathcal{ D})$ this set of functions. For $ f \in \mathcal{T}^s ( \sigma, \mu, \mathcal{ D}), $  we define 
           $$ \vert \partial_{\rho}^j f \vert_{ \sigma, \mu, \mathcal{D} } = \sup_{ x \in \mathcal{O}^s (\sigma, \mu) ; \,\rho \in \mathcal{ D}}  \max( \vert \partial_{\rho}^j f \vert, \mu \norm{ \partial_{\rho}^j \nabla_{\omega} f( x, \rho)}_s, \mu^2 \norm{ \nabla^2_{\omega_{\mathcal{ L} } \omega_{\mathcal{ L}} } \partial_{\rho}^j f( x, \rho)}),$$
        and 
        $$ [ f]^s_{\sigma, \mu, \mathcal{ D}} = \max_{j=0,1} ( \vert \partial_{\rho}^j f \vert_{ \sigma, \mu, \mathcal{D} } ).$$

   \textbf{Jet functions} For any     $ f \in \mathcal{T}^s ( \sigma, \mu, \mathcal{ D}), $  we define its jet $ f^T(x)$ as the following Taylor polynomial of $f$ at $r =0$ and $ \omega=0$ 
        	$$f^T(x)  = f( 0, \theta, 0) + d_r f(0, \theta, 0) \cdot r + d_{\omega} f( 0, \theta, 0) [\omega] + 1/2 d^2_{\omega} f(0, \theta, 0) [\omega, \omega].$$

    \textbf{Infinite matrices}
       For the elliptic variables, we denote by $\mathcal{M}_s$ the set of infinite matrices $A : \mathcal{L} \times \mathcal{L} \to \mathbb{C}$ such that $A$ maps linearly $Y_s$ into $Y_s$. We  provide $ \mathcal{M}_s$ with the operator norm
                        $$|A|_s = \norm{A}_{\mathcal{L}(Y_s,Y_s)}.$$
          We say that a matrix $A \in \mathcal{M}_s$ is in normal form if it is block diagonal
          and Hermitian, i.e.
          $$ A_{\alpha}^{\beta}= 0 \quad \text{for} \; [\alpha] \neq [\beta] \quad \text{and} \; A_{\alpha}^{\beta} = \bar{A^{\alpha}_{\beta}} \quad \text{for} \; \alpha, \beta \in \mathcal{L}.$$
          In particular, if $A \in \mathcal{M}_s$ is in normal form, its eigenvalues are
          real.

    \textbf{Normal form} A quadratic Hamiltonian function is on normal form if it
    reads
    $$ h = \Omega (\rho) \cdot r + \langle \zeta_{\mathcal{L}}, Q(\rho) \eta_{ \mathcal{L}} \rangle + 1/2 \langle \omega_{ \mathcal{ F}}, K(\rho) \omega_{ \mathcal{ F}} \rangle$$
    for some vector function $\Omega ( \rho) \in \mathbb{R}^n$, some matrix functions $ Q (\rho) \in \mathcal{M}_s$ on
    normal form and $ K ( \rho)$ is a matrix $ \mathcal{ F} \times \mathcal{F} \to \mathbb{C} $ symmetric
    in the following sense:$ K_{\alpha}^{\beta} = \, ^tK^{\alpha}_{\beta}.$

    \textbf{Poisson brackets} The Poisson brackets of two Hamiltonian functions is defined by 
    	$$ \{ f,g\} = \nabla_{ \theta} f \cdot \nabla_r g - \nabla_{r} f \cdot \nabla_{\theta}  g - i \langle \nabla_{ \omega} f , J \nabla_{\omega} g \rangle. $$
    \begin{rem}
    	A function $f$ is preserved under the flow $u(t)$ if and only if it commutes with $h$ i.e. $ \{ f,h\} =0.$ By this, we have 
    	                   $$ \{ \mathbb{L}, h \} = \{\mathbb{M},h \} =0.$$ 
    	
    \end{rem}

   \textbf{Hypothesis A0} There exists a constant $ C>0$ such that 
           $$ \vert \Lambda_a  -\vert w_a \vert^2 \vert \le C, \, \forall a \in \mathcal{L}.$$

   \textbf{Hypothesis A1} 
      \begin{align*}
      	 \vert \Lambda_{a} \vert  &\ge \delta, \quad \forall a \in \mathcal{L};\\
      	  \vert \Im \Lambda_{a} \vert  &\ge \delta, \quad \forall a \in \mathcal{ F};\\
      	 \vert \Lambda_{a}  - \Lambda_{b}\vert  &\ge \delta, \quad \forall a,b \in \mathcal{ Z}, \; [a] \neq [b];\\ 
      	 \vert \Lambda_{a}  +\Lambda_{b}\vert  &\ge \delta, \quad \forall a,b \in \mathcal{ L}.
      \end{align*}

  \textbf{Hypothesis A2} There exists $ \delta >0$ such that for all $ \Omega$ $ \delta \textendash$close to $ \Omega_0$ in $ C^1$ norm and for all $ k \in \mathbb{Z}^n \backslash \{0 \}:$ 
       \begin{enumerate}
       	\item either $$ \vert \Omega ( \rho) \cdot k \vert \ge \delta \quad \forall \rho \in \mathcal{D},$$
       	or there exists a unit vector $ z = z(k) \in \mathbb{R}^n$ such that 
       	    $$ \left( \nabla_{\rho} \cdot z \right) ( \Omega ( \rho) \cdot k) \ge \delta \quad \forall \rho \in \mathcal{ D};$$
       	 \item for all $ a \in \mathcal{ L}$
       	either $$ \vert \Omega ( \rho) \cdot k + \Lambda_{a} \vert \ge \delta \quad \forall \rho \in \mathcal{ D},$$
       	or there exists a unit vector $ z = z(k) \in \mathbb{R}^n$ such that 
       	$$ \left( \nabla_{\rho} \cdot z \right) ( \Omega ( \rho) \cdot k + \Lambda_{a}) \ge \delta \quad \forall \rho \in \mathcal{ D};$$
       	\item for all $ \alpha, \, \beta \in \mathcal{L}$ and $ a \in [ \alpha], \, b \in [ \beta]$ 
       	either $$ \vert \Omega ( \rho) \cdot k  + \Lambda_{a} \pm \Lambda_{b}\vert \ge \delta \quad \forall \rho \in \mathcal{ D},$$
       	or there exists a unit vector $ z = z(k) \in \mathbb{R}^n$ such that 
       	$$ \left( \nabla_{\rho} \cdot z \right) ( \Omega ( \rho) \cdot k + \Lambda_{a} \pm \Lambda_{b}) \ge \delta \quad \forall \rho \in \mathcal{ D};$$
       	\item for all $ a, b \in \mathcal{ F} $
       	   $$ \vert \Omega ( \rho) \cdot k  + \Lambda_{a} \pm \Lambda_{b}\vert \ge \delta.$$
       \end{enumerate}
   
  \begin{thm}[KAM theorem]\label{mainthm}
  	Assume that hypothesis A0, A1, A2 are satisfied, $f \in \mathcal{T}^s( \sigma, \mu, \mathcal{ D}),$ $f$ commutes with $\mathcal{L}, \, \mathcal{M}$ and $ s>1/2.$ Let $ \gamma >0,$ there exists a constant $C_0$ such that if 
  	   \begin{align}
  	   	    [ f]^s_{ \sigma, \mu, \mathcal{ D}} \le C_0 \delta, \quad  \varepsilon :=  [ f^T]^s_{ \sigma, \mu, \mathcal{ D}} \le C_0 \delta^{1+ \gamma},
  	   	\end{align} 
    then there exists a Cantor set $ \mathcal{ D'} \subset \mathcal{ D}$ asymptotically of full measure \textup{(}i.e. $meas ( \mathcal{ D} \setminus \mathcal{ D'}) \rightarrow 0$ when $ \varepsilon \rightarrow 0$\textup{)} and there exists a symplectic change of variables $ \Phi : \, \mathcal{O}^s( \sigma/2, \mu/2) \to \mathcal{O}^s(\sigma, \mu) $ such that for all $ \rho \in \mathcal{ D'}$
           $$ ( h_0 +f ) \circ \Phi = \tilde{h} +g$$
     with $ \tilde{h} = \Omega(\rho) \cdot r + \langle \zeta_{\mathcal{ L}}, Q(\rho) \eta_{\mathcal{ L}} \rangle +1/2 \langle \omega_{\mathcal{ F}}, K (\rho) \omega_{ \mathcal{ F}} \rangle$ on normal form, and $ g \in \mathcal{T}^s( \sigma/2, \mu/2, \mathcal{ D'}) $ with $ g^T \equiv 0.$ Furthermore there exists $C>0$ such that for all $ \rho \in \mathcal{ D'}$
      $$ \vert \Omega - \Omega_0 \vert \le C \varepsilon, \quad \vert Q - diag \left( \Lambda_{a}, \, a \in \mathcal{ L} \right) \vert \le C \varepsilon, \quad \vert JK - diag \left( \Lambda_{a}, a \in \mathcal{ F}\right) \vert \le C \varepsilon.$$
     As a dynamic consequence $ \Phi \left( \{0\} \times \mathbb{T}^n \times \{ 0\}\right)$ is an invariant torus for $h_0 +f$and this torus is linearly stable if and only if $ \mathcal{ F} = \varnothing$  \textup{(}see  \cite{maintheorem} \textup{)}. 
  \end{thm}
Here, the matrix $J$ is of the form,  
$$ \begin{pmatrix}
0 & -I\\
I &0
\end{pmatrix}$$
where $I$ is identity matrix of size $ \#F$.
\begin{rem}
   In \cite{maintheorem}, they constrained $f$ in a restricted class instead of using commutation of $f$ with $\mathcal{L}, \, \mathcal{M}$ since they considered a system of coupled NLS equation with more complicated nonlinearities.  
\end{rem}

\section{Applications}

 \textbf{The Birkhoff normal form procedure.} We recall a result proved in \cite{quinticNLSeq}.\\
 
 \begin{prop}\label{Birkhoff} 
 	There exist a canonical change of variable $\tau$ from $ \mathcal{O}^s( \sigma, \mu)$ into $ \mathcal{O}^s( 2\sigma, 2\mu)$ such that  
 	           $$ \bar{h} = h \circ \tau = N + Z_6 +R_{10}, $$
 	   where 
 	   \begin{itemize}
 	   	\item $ N$ is the term $ N ( I) = \sum_{ j \in \mathbb{Z} } j^2 |a_j|^2;$
 	   	\item $Z_6$ is the homogeneous polynomial of degree 6
 	   	             \begin{equation*}
 	   	            Z_6 = \sum_{\mathcal{R}} a_{j_1}a_{j_2}a_{j_3}b_{\ell_1}b_{\ell_2}b_{\ell_3}
 	   	             \end{equation*} 
 	   	         where \\
 	   	         $ \mathcal{R} = \{ (j, \ell ) \in \mathbb{Z}^3 \times \mathbb{Z}^3 \, s.t\, j_1 +j_2 + j_3 = \ell_1 + \ell_2 +\ell_3, \quad j_1^2 +j_2^2 + j_3^2 = \ell_1^2 +\ell_2^2 +\ell_3^2 \};$ 
 	   	\item $R_{10}$ is the remainder of order 10, i.e a Hamiltonian satisfying 
 	   	          $$ \norm{ X_{R_{10}} (x)}_s \le C \norm{ x}^9_s$$
 	   	       for all $ x \in \mathcal{O}^s(\sigma, \mu);$
 	   	\item $ \tau$ is close to the identity: there exists a constant $C$ such that 
 	   	       $$  \norm{ \tau (x) -x} \le C \norm{x}^2$$
 	   	      for all $ x \in \mathcal{O}^s(\sigma, \mu).$
 	   	            
 	   \end{itemize}
\end{prop}
Henceforth, since we do not care about constant, we shall write $a \lesssim b $ in order to say $ a \leqslant C b.$

\textbf{Persistence of 2 dimensional tori.} 

Firstly, we want to study the persistence of the two dimensional invariant torus $\mathbf{T}^2_{\nu \rho}(p,q)$ for equation \eqref{intro1} for $\nu $ small.  
Choose \\
$\begin{cases}
	a_p &= \left( \nu \rho_1 +r_1(t) \right)^\frac{1}{2} e^{i \theta_1 (t)} =: \sqrt{I_p}e^{i \theta_1(t)} \\
	a_q &= \left( \nu \rho_2 +r_2(t) \right)^\frac{1}{2} e^{i \theta_2 (t)} =: \sqrt{I_q}e^{i \theta_2(t)}\\
	a_j &= \zeta_j \hspace{1cm} j \neq p,q,
\end{cases}$\\
 where $  \{\rho_1 , \rho_2\} \in [1,2]^2 = \mathcal{ D} $ and $ \nu$ is a small parameter.
 The canonical symplectic structure now becomes 
 $$ -i d\zeta \wedge d \eta - dI \wedge d\theta $$
 with $I= ( I_1, I_2), \, \theta= ( \theta_1, \theta_2) ,\, \zeta= (\zeta_j)_j$ and $ \eta = (\eta_j)_j = ( \bar{ \zeta}_j)_j.$ \\
 Let 
 $$ \mathbf{T}^{lin}_{\rho} := \{ (I, \theta, \zeta) | \vert I -\nu \rho \vert =0, \, \vert \Im \theta \vert < \sigma , \, \norm{ \zeta}_s =0 \}$$
 and its neighborhood 
 $$ \mathbf{T}_{\rho} ( \nu, \sigma, \mu, s) := \{ (I, \theta, \zeta) | \vert I -\nu \rho \vert < \nu \mu^2, \, \vert \Im \theta \vert < \sigma , \, \norm{ \zeta}_s < \nu^{1/2}\mu \}.$$
We want to study the persistence of torus $\mathbf{T}_{\rho} ( \nu, \sigma, \mu, s)$. Indeed we have 
$$ \mathbf{T}_{\rho}  ( \nu, \sigma, \mu, s) \approx \mathcal{O}^s(\sigma, \nu^{1/2}\mu) = \{ (r, \theta, \zeta) | \vert r \vert < \nu \mu^2, \, \vert \Im \theta \vert < \sigma, \, \norm{ \zeta}_s < \nu^{1/2} \mu \}.$$
By Theorem \ref{Birkhoff}  we have 
     $$ h \circ \tau = N + Z_6 +R_{10}.$$
We see that the term $N$ contributes to the effective part and the term $R_{10}$ contributes to the remainder term $f.$ So we just need to focus on  the term $Z_6.$ Let us split it: 
     $$ Z_6 = Z_{0,6} + Z_{1, 6} + Z_{2,6} +Z_{3,6}.$$
  Here $ Z_{0,6}, \, Z_{1, 6}, \, Z_{2,6} $ are homogeneous polynomial of degree 6 which contains respectively external modes  of order $0,1,2.$ $Z_{3,6}$ is an homogeneous polynomial of degree 6 contains external modes of at least order 3,this term contributes the remainder term. 
  \\
  Thank to Lemma 2.2 on \cite{quinticNLSeq}, the term  $ Z_{1,6} =0.$ We have 
        \begin{align*}
        	Z_{0,6} &=  |a_p|^6  + | a_q|^6  + 9 \left( | a_p |^4 | a_q |^2 + | a_p|^2 |a_q |^4 \right)   \\
        	        &=  ( \nu \rho_1 + r_1)^3 + ( \nu \rho_2 + r_2)^3 + 9 \left( \nu \rho_1 + r_1 \right) \left( \nu \rho_2 + r_2 \right) \left( \nu \rho_1 + r_1 +\nu \rho_2 + r_2 \right)\\
        	    &= \nu^3 (\rho_1^3 + \rho_2^3 + 9 \rho_1^2 \rho_2 + 9\rho_2^2 \rho_1) + 3\nu^2 \left(r_1( \rho_1^2 + 6 \rho_1 \rho_2 + 3\rho_2^2 ) +r_2( \rho_2^2 + 6 \rho_1 \rho_2 + 3\rho_1^2 ) \right)  \\&+ \text{jet free}
        \end{align*}
     where the notation "jet free" means that the remaining Hamiltonian has a vanishing jet. 
 For the term $ Z_{2,6},$ there are two cases that can happen.\\ 
  \textbf{First case}\\ We assume that there is no solution\footnote{it happens when q-p is odd} $ \{s, t\} \neq \{p, q\}$ for 
          \begin{equation}
           \begin{cases} \label{st}
             2p + s  &= 2q +t\\
             2p^2 +s^2 &= 2q^2 +t^2.
           \end{cases}           
           \end{equation}
   Hence 
             $$ Z_{2,6} = Z_{2,6}^1= 9 \left(| a_p |^4 + |a_q|^4 + 4 |a_p|^2 | a_q|^2 \right)\sum_{j \neq p,q} | a_j |^2 = 9 \nu^2 \left( \rho_1^2 + \rho_2^2 + 4 \rho_1 \rho_2 \right) \sum_{j \neq p,q} | \zeta_j |^2 + \text{jet free}.$$
Hence
        $$ h \circ \tau = h^e + R$$
 where the effective Hamiltonian $h^e$ reads
\begin{align*}
	h^e &=\left( p^2 + 3 \nu^2 \left( \rho_1^2 + 3 \rho_2^2 + 6 \rho_1 \rho_2\right)\right) r_1 + \left( q^2 + 3 \nu^2 \left( \rho_2^2  + 3 \rho_1^2 + 6 \rho_1 \rho_2 \right) \right) r_2 \\&+  \sum_{j} \left( j^2 + 9\nu^2 \left( \rho_1^2 + \rho_2^2 + 4\rho_1 \rho_2 \right) \right)| \zeta_j |^2\\ 
	&=\Omega(\rho) \cdot r + \sum_{j \neq p,q} \Lambda_j | \zeta_j |^2 
\end{align*}	 
where 
     $$ \Omega ( \rho) = 
       \begin{pmatrix}
       p^2 + 3 \nu^2 \left( \rho_1^2 + 3 \rho_2^2 + 6 \rho_1 \rho_2\right)\\
       q^2 + 3 \nu^2 \left( \rho_2^2  + 3 \rho_1^2 + 6 \rho_1 \rho_2 \right)
     \end{pmatrix}$$
and 
         $$ \Lambda_{j} = j^2 + 9\nu^2 \left( \rho_1^2 + \rho_2^2 + 4\rho_1 \rho_2 \right).$$
 The remainder term $ R$ reads 
    \begin{align*}
    	R&= R_{10} + Z_{3,6}  + 3 \nu \rho_1 r_1^2 + r_1^3 + 3 \nu \rho_2 r_2^2 + r_2^3 + 9 r_1 r_2 ( r_1 +r_2) \\ &+ \left( r_1^2 + r_2^2 + 2\nu( \rho_1 +2\rho_2 ) r_1 + 2\nu( \rho_2 +2 \rho_1 ) r_2 \right) \sum_{j \neq p,q} | \zeta_j|^2 .     
    \end{align*}	
   In order to work on $ \mathcal{ O}^s ( \sigma, \mu) $ we use the rescaling 
               \begin{equation}
                \Psi: r \mapsto  \nu r, \; \zeta \mapsto \nu^{1/2} \zeta \label{rescaling}.
                \end{equation} 
   The symplectic structure now becomes 
                             $$ -\nu dr \wedge d\theta -i \nu d \zeta \wedge d \eta.$$
        By definition, this change of variables send $ \mathcal{ O}^s ( \sigma, \mu)$ to a neighborhood of $ \mathbf{T}_{\rho}(\nu,\sigma,\mu,s).$ Since $\tau $ is close to identity, the change of variables $ \Phi_{\rho}=\tau \circ \Psi$ sends $ \mathcal{ O}^s ( \sigma, \mu)$ to  $ \mathbf{T}_{\rho}(\nu,2\sigma,2\mu,s).$  By this change of variables, we have
                         $$ h \circ \Phi_{\rho}-C= ( h^e + R) \circ \Psi =  \nu h_0 + \nu f$$
        where $C$ is a constant, $h_0$ and $f$ are defined by
                          $$ h_0 = \frac{1}{\nu} h^e \circ \Psi \quad \; f = \frac{1}{\nu} R \circ \Psi.$$                                         
   By Theorem \ref{Birkhoff}, $R_{10} \in \mathcal{T}^s(\sigma, \nu^{1/2}\mu, \mathcal{D}).$ We check that the rest part of $f$ is in $ \mathcal{T}^s(\sigma,\mu, \mathcal{D}).$ By construction, $f$ commutes\footnote{since h commutes with $ \mathbb{L}$, $ \mathbb{M}$ and all the changes of variables are symplectic} with $ \mathbb{L}$ and $ \mathbb{M}$. 
   For estimating the norm of $f,$ notice that $R$ contains only term of order at least 3 in $\nu$ and $ R^T= R^T_{10}$ is of order $9/2$ in $\nu,$ so that  
                     $$ [f]^s_{ \sigma, \mu, \mathcal{D} } \lesssim \nu^{2}$$
        and 
                     $$ [f^T]^s_{ \sigma, \mu, \mathcal{D} }  \lesssim \nu^{7/2}. $$
        
So we have proved:
\begin{thm} \label{case1}
       Assume that for $ p,q \in \mathbb{Z}$ there do not exist $s,t$ solving the equation \eqref{st}. Then, the change of variables $ \Phi_{\rho} = \tau \circ \Psi $ is real holomorphic, symplectic and analytically depending on $ \rho$ satisfying 
      \begin{itemize}
      	\item $ \Phi_{\rho} : \mathcal{ O}^s ( \sigma, \mu) \to \mathbf{T}_{\rho} ( \nu, 2\sigma, 2\nu, s) ;$
      	\item $ \Phi_{\rho}$ puts the Hamiltonian $h$ in normal form in the
      	following sense: $$ \frac{1}{\nu} ( h \circ \Phi_{\rho} - C) = h_0 + f$$
      	 where $C$ is a constant and the effective part $ h_0$ of the Hamiltonian reads 
      	 \begin{align*}
      	 &h_0= \Omega(\rho) \cdot r + \sum_{j \neq p,q} \Lambda_j | \zeta_j |^2 
      	 \end{align*}	 
      	 with 
      	 $$ \Omega ( \rho) = 
      	 \begin{pmatrix}
      	 p^2 + 3 \nu^2 \left( \rho_1^2 + 3 \rho_2^2 + 6 \rho_1 \rho_2\right)\\
      	 q^2 + 3 \nu^2 \left( \rho_2^2  + 3 \rho_1^2 + 6 \rho_1 \rho_2 \right)
      	 \end{pmatrix}$$
      	 and 
      	 $$ \Lambda_{j} = j^2 + 9\nu^2 \left( \rho_1^2 + \rho_2^2 + 4\rho_1 \rho_2 \right);$$
      	  \item The remainder term $f$ belongs to $ \mathcal{T}^s ( \sigma, \mu, \mathcal{ D})$   and satisfies 
      	        $$ [f]^s_{ \sigma, \mu, \mathcal{D} } \lesssim \nu^{2}$$
      	        and 
      	        $$ [f^T]^s_{ \sigma, \mu, \mathcal{D} }  \lesssim \nu^{7/2}. $$
      \end{itemize}

\end{thm}
\textbf{Second case}\\
  Assume that there are\footnote{in this case, $\{p,q,s,t\}$ is of the form $\{p, p+2n,p+3n,p-n\}$ } $ s, t \neq p, q$ solving \eqref{st}, hence
         \begin{align*}
        Z_{2,6} = Z_{2,6}^1  + 9 (a_p^2 a_{s} b_q^2 b_{t}  + b_p^2 b_{s} a_q^2 a_{t}) = Z_{2,6}^1 + Z_{s,t}
        \end{align*}
   For the second term, let us rewrite it
       $$  9( \nu \rho_1 +r_1) ( \nu \rho_2 +r_2) \left( e^{ 2i ( \theta_1 -\theta_2 ) } \zeta_{s} \eta_{t} +  e^{ -2i ( \theta_1 -\theta_2 ) } \eta_{s} \zeta_{t} \right)$$
   The effective part of this term is just given by 
       $$ 9\nu^2 \rho_1 \rho_2\left( e^{ 2i ( \theta_1 -\theta_2 ) } \zeta_{s} \eta_{t} +  e^{ -2i ( \theta_1 -\theta_2 ) } \eta_{s} \zeta_{t} \right) .$$
   Notice that 
       $$ \{I_s, \zeta_s \eta_t + \eta_s \zeta_t \} = \{ I_t, \zeta_s \eta_t + \eta_s \zeta_t\} = 0.$$
   This gives us a clue that the above term does not effect to the stability of the solution. \\
   In order to kill the angles, we introduce the symplectic change of variables $ \Psi_{angles} : \mathcal{O}^s (\sigma, \mu) \to \mathcal{O}^s (\sigma, \mu),$ $( r_1, r_2, \theta, \zeta) \mapsto \left( r'_1, r'_2, \theta, \zeta' \right)$ defined by 
     $$\begin{cases}
     	 \zeta_s' &= e^{2i ( \theta_1 -\theta_2)} \zeta_s \\
     	 \zeta_t' &= \zeta_t \\
     	 \zeta_j' &= \zeta_j, \quad j \neq s,t, p,q \\
     	 r_1' &= r_1 -2|\zeta_s|^2 \\
     	 r_2' &= r_2 + 2|\zeta_s|^2. 
     \end{cases}$$
     By this change of variables 
        $$ \tilde{h} = \bar{h} \circ \Psi_{angles} = C + h^e +R.$$      
    Here $C$ is a constant given by 
         $$ C= \nu^3 (\rho_1^3 + \rho_2^3 + 9 \rho_1^2 \rho_2 + 9\rho_2^2 \rho_1) + 9(\nu p^2 \rho_1 + \nu q^2 \rho_2).$$
        The effective Hamiltonian $h^e$ reads 
    \begin{align*}
    	h^e &= \left( p^2 + 3 \nu^2 \left( \rho_1^2 + 3 \rho_2^2 + 6 \rho_1 \rho_2\right)\right) r_1' + \left( q^2 + 3 \nu^2 \left( \rho_2^2  + 3 \rho_1^2 + 6 \rho_1 \rho_2 \right) \right) r_2' \\&+  \sum_{j \neq p,q,s,t} \left( j^2 + 9\nu^2 \left( \rho_1^2 + \rho_2^2 + 4\rho_1 \rho_2 \right) \right)| \zeta_j' |^2  + \left( t^2 + 9\nu^2 \left( \rho_1^2 + \rho_2^2 + 4\rho_1 \rho_2 \right) \right)| \zeta_t' |^2  \\&+ \left( s^2 +2p^2 - 2 q^2 + \nu^2 \left(21 \rho_2^2 -3 \rho_1^2  + 36\rho_1 \rho_2 \right) \right)| \zeta_s' |^2   + 9\nu^2 \rho_1 \rho_2 ( \zeta_s' \eta_t' + \eta_s' \zeta_t').
    \end{align*}	
     It is on normal form 
    \begin{align*} 
    	  \Omega(\rho) \cdot r + \sum_{j \neq p,q,s,t} \Lambda_j | \zeta_j' |^2 + \Lambda_s |\zeta_s'|^2 + \Lambda_t |\zeta_t'|^2 + 9\nu^2 \rho_1 \rho_2 ( \zeta_s' \eta_t' + \eta_s' \zeta_t')
    \end{align*}	 
    where 
    $ \Omega (\rho)$ and $\Lambda_j$ are defined as in the first case except 
          $$ \Lambda_s = t^2 + \nu^2 \left(21 \rho_2^2 -3 \rho_1^2  + 36\rho_1 \rho_2\right)  .$$
   In order to diagonalize $h^e$, we use a symplectic change of variables of the form
                   $$ \begin{cases}
                     \zeta_{t^+} &= \frac{1}{\sqrt{1+ \alpha^2}}(\zeta_t' + \alpha \zeta_s') \\
                      \zeta_{t^-} &= \frac{1}{\sqrt{1+ \alpha^2}}(\zeta_s' - \alpha \zeta_t')
                   \end{cases}$$ 
        with $ \alpha= \frac{  -2\rho_1^2 +2\rho_2^2 +\sqrt{ 4\rho_1^4 +2\rho_1^2\rho_2^2+4\rho_2^4  } }{ 3\rho_1 \rho_2}.$           
    Then $h^e$ can be diagonalized as
               $$ \Omega(\rho) \cdot r + \sum_{j \neq p,q,s,t} \Lambda_j | \zeta_j |^2 + \Lambda_{t^+} |\zeta_{t^+}|^2 + \Lambda_{t^-} |\zeta_{t^-}|^2$$
       where                  
                   $$
                   \begin{cases}                  
                   \Lambda_{t^+} &= \Lambda_t - 9\nu^2 \rho_1\rho_2 \alpha\\
                   \Lambda_{t^-} &= \Lambda_s +9\nu^2 \rho_1 \rho_2 \alpha.
                   \end{cases}
                   $$
   The remainder term $ R$ reads 
   \begin{align*}
   	R&= R_{10} \circ \Psi_{angles} + Z_{3,6} \circ \Psi_{angles}+ 3 \nu \rho_1 r_1^2 + r_1^3 + 3 \nu \rho_2 r_2^2 + r_2^3  \\&+9 r_1 r_2 ( r_1 +r_2) + \left( r_1^2 + r_2^2 + 2\nu( \rho_1 +2
   	\rho_2 ) r_1 + 2\nu( \rho_2 +2 \rho_1 ) r_2 \right) \sum_{j \neq p,q} | \zeta_j|^2 
   \end{align*}	
   with $ r_1 = r_1' +2|\zeta_s|^2, \; r_2 = r_2' -2 |\zeta_s|^2.$ \\
   Using the rescaling $\Psi$ introduced in \eqref{rescaling}, we get
                         $$ ( h^e + R) \circ \Psi = \nu h_0 + \nu f.$$
    Since $\Psi_{angles}: \mathcal{O}^s ( \sigma,  \mu) \to \mathcal{O}^s ( \sigma,  3\mu)$ and $\tau$ is closed to identity, we have $\tau \circ \Psi_{angles} \circ \Psi : \mathcal{O}^s ( \sigma,\mu) \to \mathbf{T}_{\rho} ( \nu, 2\sigma, 4\mu, s). $
   The study of $f$ is the same as in the previous case. Then we get: 
\begin{thm} \label{case2}
	Assume that $ p,\,q, \,s, \, t$  satisfy the equation \ref{st}. The change of variables $ \Phi_{\rho} = \tau \circ \Psi_{angles} \circ \Psi  $ is a real holomorphic transformations, analytically depending on $ \rho$ satisfying 
	\begin{itemize}
		\item $ \Phi_{\rho} : \mathcal{ O}^s ( \sigma,\mu) \to \mathbf{T}_{\rho} ( \nu, 2\sigma, 4\mu, s) ;$
		\item $ \Phi_{\rho}$ puts the Hamiltonian $h$ in normal form in the
		following sense: $$ \frac{1}{\nu} ( h \circ \Phi_{\rho} - C) = h_0 + f$$
		where $C$ is a constant and the effective part $ h_0$ of the Hamiltonian reads 
		\begin{align*}
		&h_0= \Omega(\rho) \cdot r + \sum_{j \neq p,q, s,t} \Lambda_j | \zeta_j |^2  + \Lambda_{t^+} |\zeta_{t^+}|^2 + \Lambda_{t^-} |\zeta_{t^-}|^2
		\end{align*}	 
		with
		$$ \Omega ( \rho) = 
		\begin{pmatrix}
		p^2 + 3 \nu^2 \left( \rho_1^2 + 3 \rho_2^2 + 6 \rho_1 \rho_2\right)\\
		q^2 + 3 \nu^2 \left( \rho_2^2  + 3 \rho_1^2 + 6 \rho_1 \rho_2 \right)
		\end{pmatrix}$$
		and 
		$$ \Lambda_{j} = j^2 + 9\nu^2 \left( \rho_1^2 + \rho_2^2 + 4\rho_1 \rho_2 \right),$$
		
		\item The remainder term $f$ belongs to $ \mathcal{T}^s ( 1, 1, \mathcal{ D})$   and satisfies 
		$$ [f]^s_{ \sigma, \mu, \mathcal{D} } \lesssim \nu^{2}$$
		and 
		$$ [f^T]^s_{ \sigma, \mu, \mathcal{D} }  \lesssim \nu^{7/2}. $$
	\end{itemize}

\end{thm}
Now we can finish the proof of Theorem \ref{thm2}.

\textit{Proof of Theorem \ref{thm2}}. By Theorem \ref{case1} and \ref{case2}, there exists a symplectic change of variables $\Phi_{\rho}$, on a asymtotical set $\mathcal{D}_{\nu} \mathcal{D} = [1,2]^2,$ puts the Hamiltonian $h= N+P$ in normal form $h_0 + f,$ that satisfies,(see the appendix A) the hypotheses of KAM theorem \ref{mainthm} for $ \delta = \nu^2$, $ \varepsilon = \nu^{7/2}= \delta^{7/4}$ and $ \Omega_0=\omega = (p^2, q^2) + O(\nu^2).$ So by KAM theorem, since the hyperbolic set $\mathcal{F}$ is empty, the torus\footnote{here we choose $\sigma=1$} 
$$ \mathbf{T}^{lin}_{\rho} := \{ (I, \theta, \zeta) | \vert I -\nu \rho \vert =0, \, \vert \Im \theta \vert < 1 , \, \norm{ \zeta}_s =0 \}$$
 is linear stable. Here we denote $ I=(I_p, I_q)$.
          \begin{flushright}
           	$\Box$
           \end{flushright}

\textbf{Persistence of 3 dimensional tori.}
Assume that \\
$$\begin{cases}
a_p &= \left( \nu \rho_1 +r_1(t) \right)^\frac{1}{2} e^{i \theta_1 (t)} =: \sqrt{I_p}e^{i \theta_1(t)} \\
a_q &= \left( \nu \rho_2 +r_2(t) \right)^\frac{1}{2} e^{i \theta_2 (t)} =: \sqrt{I_q}e^{i \theta_2(t)}\\
a_m &= \left( \nu \rho_3 +r_3(t) \right)^\frac{1}{2} e^{i \theta_3 (t)} =: \sqrt{I_m}e^{i \theta_3(t)}\\
a_j &= \zeta_j \qquad j \in \mathbb{Z} \setminus \{ p,q,m\}
\end{cases}$$\\
where $  \rho = (\rho_1 , \rho_2, \rho_3) \in \mathcal{D} \subset \mathbb{R}^3 $ and $ \nu$ is a small parameter. The canonical symplectic structure now becomes 
$$ -i d\zeta \wedge d \eta - dI \wedge d\theta $$
with $I= ( I_p, I_q, I_m), \, \theta= ( \theta_1, \theta_2 , \theta_3), \, \zeta= (\zeta_j)_{ j \in \mathbb{Z} \setminus \{ p,q,m\}}$ and $ \eta = (\eta_j)_{ j \in \mathbb{Z} \setminus \{ p,q,m\}} = ( \bar{ \zeta}_j)_{ j \in \mathbb{Z} \setminus \{ p,q,m\}}.$ \\
The same as in two-modes case, we have
$$ \bar{h} := h \circ \tau= N + Z_6 +R_{10}.$$
We see that as in the previous case, the term $N$ contributes to the effective Hamiltonian $h_0$ and the term $R_{10}$ contributes to the remainder term $f.$ So we just need to focus on  the term $Z_6.$ Let us split it: 
$$ Z_6 = Z_{0,6} + Z_{1, 6} + Z_{2,6} +Z_{3,6}.$$
Here, $Z_{0,6}$ is homogeneous polynomial of degree 6 which just contains inner modes $ (p,q,m)$;  $Z_{1,6}$, $Z_{2,6}$ are homogeneous polynomials of degree 6 which contain outer modes of order $1$ and $2$. $Z_{3,6}$ is an homogeneous polynomial of degree 6 contains outer modes of at least order $3,$ this term contributes the remainder term. 
We have: 
$$	Z_{0,6} =  |a_p|^6  + | a_q|^6  + | a_m|^6 + 9  \sum_{ j, \ell  \in \{ p,q,m\} }|a_j|^4 |a_{\ell}|^2   + 36 |a_p|^2 |a_q|^2 |a_m|^2  $$ 
Even if it looks a bit more complicated, we deal with $Z_{0,6}$ as in the previous case. 
We assume that there is no solution to \eqref{introell}, so that $Z_{1,6}=0.$
For $Z_{2,6},$ we have
\begin{align*} Z_{2,6} &= \sum_{j_1, j_2, \ell} |a_{j_1}|^2 |a_{j_2}|^2 |a_{\ell}|^2 + \sum_{s_1, t_1 \in \mathcal{A}} \left( a_{j_3}^2 a_{s_1} b_{j_4}^2 b_{t_1}  +  b_{j_3}^2 b_{s_1} a_{j_4}^2 a_{t_1} \right) \\ &\hspace{1cm}+ \sum_{ s_2, t_2 \in \mathcal{B}} \left( a_{j_5}^2 a_{j_6} b_{j_7} b_{s_2} b_{t_2}  +  b_{j_5}^2 b_{j_6} a_{j_7} a_{s_2} a_{t_2} \right)  
	\\ &\hspace{1cm} + \sum_{ s_3, t_3 \in \mathcal{C}} \left( a_{j_9}^2 a_{s_3} b_{j_8} b_{j_{10}} b_{t_3}  +  b_{j_9}^2 b_{s_3} a_{j_8} a_{j_{10}} a_{t_3} \right) \\
	&\hspace{1cm}+ \sum_{ s_4 \in \mathcal{E}} \left( a_{j_{11}}^2 a_{j_{12}} b_{j_{13}} b_{s_4}^2  +  b_{j_{11}}^2 b_{j_{12}} a_{j_{13}} a_{s_4}^2 \right)      
\end{align*}
with $ j_i \in \{ p,q,m\},$ $s_i,t_i \notin \{p,q,m\}$ and $s_i \neq t_i.$ The sets $ \mathcal{ A}, \,\mathcal{ B},\, \mathcal{ C},\, \mathcal{ E}$ are given by
\begin{align*}
 \mathcal{A}\leftrightarrow	&\begin{cases}
		2j_3 +s_1 &= 2 j_4 + t_1\\
		2j_3^2 +s_1^2 &= 2 j_4^2 + t_1^2
	\end{cases}                         
	\qquad &
\mathcal{B}\leftrightarrow	&\begin{cases}
		2j_5 +j_6  &= j_7 + s_2+ t_2\\
		2j_5^2 +j_6^2  &= j_7^2 + s_2^2+ t_2^2
	\end{cases}\\
\mathcal{C}\leftrightarrow	&\begin{cases}
		2j_9 +s_3  &= j_8 + j_{10}+ t_3\\
		2j_9^2 +s_3^2  &= j_8^2 + j_{10}^2+ t_3^2
	\end{cases}
	\qquad &
\mathcal{E}\leftrightarrow	&\begin{cases}
		2j_{11} +j_{12}  &= j_{13} + 2s_4\\
		2j_{11}^2 +j_{12}^2  &= j_{13}^2 + 2s_4^2.
	\end{cases}
\end{align*}
Assume that $ \mathcal{A}, \mathcal{B}, \mathcal{C}, \mathcal{ E}$ are disjoint\footnote{this is the case for the example considered in theorem \ref{thm3}} i.e. there is no s or t appearing in two of these sets. We shall deal with each term one by one (in case it's not empty). 

The first term just depends on the actions, and we have
$$ |a_{j_1}|^2 |a_{j_2}|^2 |a_{\ell}|^2 = \nu^2 \rho_{j_1} \rho_{ j_2}|\zeta_{\ell}|^2 + \text{jet free}. $$ 
              
The second and the fourth term are similar, since their effective parts are all of the form 
$$ 9e^{ i \alpha } \zeta_s \eta_t + 9e^{- i \alpha} \eta_s \zeta_t.$$
The idea to deal with these two terms is the same as that in the two-modes case. Since  $$\{ I_s +I_t, \zeta_s \eta_t \} = \{ I_s + I_t, \zeta_t \eta_s\} =0,$$
these terms do not affect the stability of the flow. Since $ \mathcal{A}, \mathcal{B}, \mathcal{C}, \mathcal{ E}$ are disjoint, and as in the two-modes case, a change of variables that used to deal with a pair ${s,t}$ only affect that modes, i.e the changes of variables commute.  We call $ \Phi_1$ the composition of all changes of variables used to deal with the sets $\mathcal{ A}$ and $\mathcal{ C}$. 

For the third term,  its effective parts are of the form 
$$ 18\nu^2 \rho_{j_5} \sqrt{ \rho_{ j_6} \rho_{ j_7}}(e^{ i \alpha } \zeta_s \zeta_t + e^{- i \alpha} \eta_s \eta_t)$$
where $ \alpha = \theta_{j_7} - \theta_{j_6} -2\theta_{j_5}.$ For explicitness, we will consider the case $ j_5 = p, \; j_6= q, \; j_7=m,$ and $s,\,t$ solve the following equation 
\begin{equation}
\label{zst3} \begin{cases}
2p + q &= m +s + t \\
2p^2 + q^2 &=m^2 +s^2 +t^2,
\end{cases}
\end{equation}  
then $ \alpha = \theta_{3} - \theta_{2} -2\theta_{1}.$ An example for this could be $ ( p,q,m, s, t) = ( -3, 10, -6, 1, 9).$
In order to kill the angles, we introduce the symplectic change of variables $\Psi_{ang,1}:\mathcal{ O}^s(\sigma, \mu) \to \mathcal{ O}^s(\sigma, 3\mu);\,( r, \theta, \zeta) \mapsto ( r', \theta, \zeta')$
defined by \\
$$\begin{cases}
\zeta_s' &= ie^{-i\alpha} \eta_s \qquad \eta_s' = ie^{i\alpha} \zeta_s \\
\zeta_t' &= \zeta_t  \hspace{1.8cm} \eta_t' = \eta_t \\
\zeta_j' &= \zeta_j, \hspace{1.55cm}\eta_j' = \eta_j \quad j \neq s,t, p,q \\
r_1' &= r_1 +2|\zeta_s|^2 \\
r_2' &= r_2 + |\zeta_s|^2,        \\
r_3' &=r_3 - |\zeta_s|^2.        
\end{cases}$$
The effective part related to $ s, t$ is of the form 
\begin{equation}
\Lambda_s | \zeta_s'|^2 + \Lambda_t | \zeta_t'|^2  - 18i\nu^2 \rho_1 \sqrt{ \rho_2 \rho_3 }(  \zeta_s' \eta_t' + \eta_s'\zeta_t')\label{stB}
\end{equation}
where $$ \Lambda_t= t^2 + 9\nu^2 ( \rho_1^2  + \rho_2^2 +\rho_3^2  + 4 \rho_1 \rho_2 + 4 \rho_2 \rho_3 + 4 \rho_3 \rho_1 )$$
and   $$ \Lambda_s = t^2 + 3\nu^2 ( -\rho_1^2 + \rho_2^2 + 5 \rho_3^2 -6 \rho_1 \rho_2 +12 \rho_2 \rho_3 + 6 \rho_3 \rho_1  ).$$
Denoting $ a = \frac{ \Lambda_t - \Lambda_s}{2}$ and $b = \frac{ \Lambda_t + \Lambda_s}{2},$ we diagonalize \eqref{stB} by the symplectic change of variables\footnote{$\sqrt{-1}=i$}
$$
\begin{cases}
\zeta_{t^-} &= \frac{1}{\sqrt{1- \alpha^2}} ( \zeta_s' - i \alpha \zeta_t')  \quad  \eta_{t^-} = \frac{1}{\sqrt{1- \alpha^2}} ( \eta_s' - i \alpha \eta_t')\\ 
\zeta_{t^+} &= \frac{1}{\sqrt{1- \alpha^2}}( \zeta_t' + i \alpha \zeta_s')   \quad  \eta_{t^+} = \frac{1}{\sqrt{1- \alpha^2}}( \eta_t' +i \alpha \eta_s')
\end{cases}
$$ 	 
where  $$\alpha =- \frac{a - \sqrt{a^2-18^2\nu^4 \rho_1^2 \rho_2 \rho_3}}{ \nu^2 \rho_1  \sqrt{ \rho_2 \rho_3 }} . $$
Then \eqref{stB} becomes
$$ \Lambda_{t^+}| \zeta_{t^+}|^2 +  \Lambda_{t^-}| \zeta_{t^-}|^2  $$
where $\Lambda_{t^{\pm}} = b \pm \sqrt{a^2-18^2\nu^4 \rho_1^2 \rho_2 \rho_3}.$ We see that two modes $ t^+ ,\, t^-$ correspond to hyperbolic direction if and only if $ a^2-18^2\nu^4 \rho_1^2 \rho_2 \rho_3 < 0$, a condition related to the choice of $\rho.$ Precisely, for $ \rho \in \mathcal{D}_1 = [1,2]^3$, we have $\Lambda_{t^{\pm}} \in \mathbb{R}$ while for $\rho=(2, 1, 9)$ we have $a=0$ and $a^2 -18^2\nu^4 \rho_1^2 \rho_2 \rho_3= -18^2\nu^4 \rho_1^2 \rho_2 \rho_3<0$. Hence, there exist $ \epsilon >0$(choose $\epsilon =10^{-2}$) such that for $ \rho \in \mathcal{ D}_2 = \mathcal{ D}_{\epsilon} = [2- \epsilon, 2+ \epsilon] \times [1- \epsilon,1+ \epsilon]  \times [9- \epsilon, 9+ \epsilon]$ we have $ | \Im \Lambda_{t^{ \pm}}| > \nu^2.$
We call $ \Phi_2$ the composition of changes of variables related to $ \mathcal{ B}.$

For the set $\mathcal{E}$, without loss of generality, assume that 
\begin{equation}
\label{zss3} \begin{cases}
2p + q &= m +2s  \\
2p^2 + q^2 &=m^2 +2s^2 .
\end{cases}
\end{equation} 
Then, using the symplectic change of variables $\Psi_{ang,2}:\mathcal{ O}^s(\sigma, \mu) \to \mathcal{ O}^s(\sigma, 2\mu);\,( r, \theta, \zeta) \mapsto ( r', \theta, \zeta')$
defined by \\
$$\begin{cases}
\zeta_s' &= e^{i\alpha/2} \zeta_s \quad \eta_s' = e^{-i\alpha/2} \eta_s \\
\zeta_j' &= \zeta_j,  \hspace{1cm}  \eta_j' = \eta_j \quad j \neq s, p,q \\
r_1' &= r_1 +|\zeta_s|^2 \\
r_2' &= r_2 +\frac{1}{2} |\zeta_s|^2       \\
r_3' &=r_3 - \frac{1}{2}|\zeta_s|^2.        
\end{cases}$$
The effective part related to $s$ becomes
\begin{equation} 
\Lambda_s | \zeta_s'|^2 + \nu^2 \rho_1 \sqrt{ \rho_2 \rho_3 } ( \zeta_s'^2 +\eta_s'^2) \label{ss}
\end{equation}
where $$\Lambda_s = 3 \nu^2 (2\rho_1^2  + \rho_2^2 -\rho_3^2  + 9 \rho_1 \rho_2 + 3 \rho_3 \rho_1 )$$
If $ \Lambda_s \neq 0,$ we can diagonalize \eqref{ss} into $ \frac{1- \beta^2}{ 1 + \beta^2}\Lambda_s | \frac{\zeta_s' + \beta \eta_s'}{ \sqrt{1 - \beta^2}} |^2$ with $\beta$ satisfying $ \Lambda_s\beta = (1- \beta^2) \nu^2 \rho_1 \sqrt{ \rho_2 \rho_3 },$ otherwise we rewrite it into $i\nu^2 \rho_1 \sqrt{ \rho_2 \rho_3 } ( \frac{\zeta_s' + i \eta_s'}{\sqrt{2}} \frac{\eta_s' + i \zeta_s'}{\sqrt{2}}),$ however $meas \{ \rho \in \mathbb{R}^3: \Lambda_s = 0\} =0$. We call $ \Phi_3$ the composition of all changes of variables related to $ \mathcal{ E}.$ 

By construction of $\Phi_i$ and definition of $\mathcal{O}^s(\sigma,\nu)$, the composition $ \Phi_3 \circ \Phi_2 \circ \Phi_1$ mapping $\mathcal{O}^s(\sigma,\nu)$ into $\mathcal{O}^s(\sigma,3\nu).$  
Using the rescaling $\Psi$ introduced in \eqref{rescaling}, as the previous case we get 
\begin{thm} \label{thmmode3}
	Assume that the equation \eqref{introell} with $j_1,j_2,j_3 \in \{ p,q,m\}$ has no solution in $\mathbb{Z}$ and $ \mathcal{A}, \, \mathcal{B},\, \mathcal{C,\,\mathcal{E}}$ are disjoint. The change of variables $ \Phi_{\rho} := \Psi \circ \Phi_3 \circ \Phi_2 \circ \Phi_1 \circ \tau  $ is a holomorphic, symplectic transformation, and analytically depending on $\rho \in \mathcal{ D}$, satisfying 
	\begin{itemize}
		\item $ \Phi_{\rho}: \mathcal{ O}^s(\sigma,\mu) \to \mathbf{T}_{\rho} ( \nu, 2\sigma,4\mu,s);$
		\item $ \Phi_{\rho}$ puts the Hamiltonian $h$ in normal form in the
		following sense: $$ \frac{1}{\nu} ( h \circ \Phi_{\rho} - C) = h_0 + f$$
		where $C$ is a constant and the effective part $ h_0$ of the Hamiltonian reads 
		\begin{align*}
			&h_0= \Omega(\rho) \cdot r + \sum_{a \in \mathcal{ Z}} \Lambda_a | \zeta_a |^2 
		\end{align*}	 
		where 
		$$ \Omega ( \rho) = 
		\begin{pmatrix}
		p^2 + 3 \nu^2 \left( \rho_1^2 + 3 \rho_2^2 + 3 \rho_3^2 +   6 \rho_1 \rho_2 + 6 \rho_1 \rho_3 +12 \rho_2 \rho_3\right)\\
		q^2 + 3 \nu^2 \left( \rho_2^2 + 3 \rho_1^2 + 3 \rho_3^2 +   6 \rho_1 \rho_2 + 6 \rho_2 \rho_3 +12 \rho_1 \rho_3 \right)\\
		m^2 + 3 \nu^2 \left( \rho_3^2 + 3 \rho_1^2 + 3 \rho_2^2 +   6 \rho_1 \rho_3 + 6 \rho_3 \rho_2 +12 \rho_2 \rho_1\right)
		\end{pmatrix}$$
		\item $ \mathcal{ Z}$ is the disjoint union $ \mathcal{ L} \cup \mathcal{ F};$ $ \mathcal{ L}$ corresponds to elliptic part, and $ \mathcal{ F}$ corresponds to hyperbolic part; 
		\item  the remainder term $f$ belongs to $ \mathcal{T}^s ( \sigma, \mu, \mathcal{ D})$   and satisfies 
		$$ [f]^s_{ \sigma, \mu, \mathcal{D} } \lesssim \nu^{2}$$
		and 
		$$ [f^T]^s_{ \sigma, \mu, \mathcal{D} }  \lesssim \nu^{7/2}. $$
	\end{itemize}
	
\end{thm} 

 
\textit{Proof of theorem \ref{thm3}.} By theorem \ref{thmmode3}, for $(p,q,m)=(-3,10,-6)$ and $\rho \in  \mathcal{D}_{\nu} \subset \mathcal{ D}_2$, there exists a symplectic change of variables $\Phi_1$ on $ \mathcal{D}_{\nu}$ puts the Hamiltonian $h= N+P$ in normal form $h_0 + f,$ that satisfies,(see appendix A) assumptions of KAM theorem \ref{mainthm} for $ \delta = \nu^2$, $ \varepsilon = \nu^{7/2}= \delta^{7/4}$ and $ \Omega_0 = \omega = ( 3^2, 10^2, 6^2) + \mathcal{ O}(\nu^2).$ So by KAM theorem, 
 the torus 
$$ \mathbf{T}^{lin}_{\rho} = \{ (I, \theta, \zeta) | \vert I -\nu \rho \vert =0, \, \vert \Im \theta \vert < 1 , \, \norm{ \zeta}_s =0 \}$$
 is linearly unstable. 
\begin{flushright}
	$\Box$
\end{flushright}


\section{Appendix A}
In this appendix, we will verify the hypothesis A0, A1, A2 of Theorem \ref{mainthm} for the Hamiltonian in our applications. The hypothesis A0, A1 is trivial, so we focus on A2. 
\subsection{Two-modes case}
\textbf{The first case} In this case, we have $ \mathcal{ F} = \emptyset$ and the other estimates are trivial. For the hypothesis A2, we recall that 
       $$ \Omega ( \rho) = 
      \begin{pmatrix}
      p^2 + 3 \nu^2 \left( \rho_1^2 + 3 \rho_2^2 + 6 \rho_1 \rho_2\right)\\
      q^2 + 3 \nu^2 \left( \rho_2^2  + 3 \rho_1^2 + 6 \rho_1 \rho_2 \right)
      \end{pmatrix}$$
      and 
      $$ \Lambda_{j} = j^2 + 9\nu^2 \left( \rho_1^2 + \rho_2^2 + 4\rho_1 \rho_2 \right).$$
      Let $ k = ( k_1, k_2) \in \mathbb{Z }^2/ \{ 0\}$ and $ z = z( k) = \frac{(k_2,k_1)}{ | k|},$ then we have 
                  \begin{align*}  ( \nabla_{\rho} \cdot z) ( \Omega (\rho) \cdot k) &= 6 \nu^2 \left( 3( \rho_1 +  \rho_2) k_2^2 + 3(\rho_2 + 3\rho_1) k_1^2 + 4 ( \rho_1 + \rho_2) k_1 k_2 \right) | k|^{-1} \\
                  	&\ge   \frac{6}{\sqrt{2}} \nu^2 |k|
                  	\end{align*}
         and 
         $$ (\nabla_{\rho} \cdot z) \Lambda_{j}  = 18 \nu^2 ( ( \rho_1 + 2 \rho_2) k_2 + ( \rho_2 + 2 \rho_1) k_1 ) | k|^{-1}
         .$$
         Choosing $ \delta = 4\nu^2 $, we get the hypothesis A2 $(1).$ Since $(\nabla_{\rho} \cdot z) (\Lambda_{j} - \Lambda_{\ell} ) =0 ,$ the estimate of small divisor $\Omega \cdot k + \Lambda_{j} - \Lambda_{\ell} $  is followed. To estimate the small divisors $ \Omega \cdot k + \Lambda_{j}$ and $ \Omega \cdot k + \Lambda_{j} + \Lambda_{\ell}$ we use the fact that f commute with both the mass $ \mathbb{ L}$ and momentum $\mathbb{ M}.$ We just need to control small divisors $ \Omega \cdot k + \Lambda_{j}$ and $ \Omega \cdot k + \Lambda_{j} + \Lambda_{\ell}$ whenever $ e^{i k \cdot \theta} \eta_j \in f  $ and $ e^{i k \cdot \theta} \eta_j \eta_{\ell} \in f$, respectively.
         We have for the mass and momentum: 
          $$ \mathbb{L} = \nu ( \rho_1 + \rho_2) + r_1 +r_2 + \sum_j | \zeta_j |^2$$
          and 
          $$ \mathbb{M} = \nu ( p\rho_1 + q\rho_2) + pr_1 +qr_2 + \sum_j j| \zeta_j |^2.$$	
           By conservation of $ \mathbb{L},$ we have 
                         $$ \{ e^{i k \cdot \theta} \eta_j,\mathbb{L} \} = ie^{i k \cdot \theta} \eta_j ( k_1 + k_2 +1) =0. $$
              Therefore, for A2 $(2)$ we just have to study the case $ k_1 + k_2 =-1.$ In this situation 
                          \begin{align*}
                          	( \nabla_{\rho} \cdot z) ( \Omega (\rho) \cdot k + \Lambda_{j}) &=  6 \nu^2 |k|^{-1} \left( 3(  \rho_1 +  \rho_2) k_2^2 + 3( \rho_2 +  \rho_1) k_1^2 + 4 ( \rho_1 + \rho_2) k_1 k_2 \right)  \\
                          	&\hspace{1cm}+ 6 \nu^2 |k|^{-1} \left(3( \rho_1 + 2 \rho_2) k_2 + 3( \rho_2 + 2 \rho_1) k_1 \right)\\
                          	    &= 6 \nu^2 |k|^{-1}  \left(  ( \rho_1 +  \rho_2) k_2^2 + ( \rho_2 +  \rho_1) k_1^2 + 2 ( \rho_1 + \rho_2)   \right) \\
                          	    &\hspace{1cm}+  6 \nu^2 |k|^{-1} \left( 3\rho_2 k_2 + 3\rho_1 k_1 - 3 ( \rho_1 + \rho_2) \right)\\ 
                          	    &=6 \nu^2 |k|^{-1} \left( 2( \rho_1 + \rho_2) k_1^2 + ( 5 \rho_1 - \rho_2) k_1 - 3  \rho_2 \right).
                          	\end{align*}
          This term is greater than $ \delta$ except the cases $ k = ( -1,0)$ and $ (0,-1).$ The conservation of $\mathbb{ M}$ gives us
                      $$ \{ e^{i k \cdot \theta} \eta_j,\mathbb{M}\} = ie^{i k \cdot \theta} \eta_j ( pk_1 + qk_2 +j)=0.$$ 
             For $ k \in \{ (-1,0), \, (0,-1)\},$ this implies $ j \in \{ p,q\},$ which is excluded. \\
             We consider the small divisor $ \Omega \cdot k + \Lambda_{j} + \Lambda_{\ell}$  in the same way. The conservation of the mass $ \mathbb{ L}$ gives us $ k_1 + k_2 = -2$ and then by computation we get $ k \in \{ ( 0,-2), ( -2, 0), ( -1, -1), (-3,1), (1,-3) \}.$ The conservation of the momentum gives us $ p k_1 + qk_2 + j + \ell =0.$ We have 
                        $$ \Omega \cdot k + \Lambda_{j} + \Lambda_{\ell} = N ( p,q, j, \ell) + \mu (\rho, k, )$$
                  where $ N ( p, q, j, \ell) = p^2 k_1 + q^2 k_2 + j^2 + \ell^2$ and $\mu (\rho)$ very small for $ |k| \le 4.  $
                  We see that $ N ( p, q, j, \ell) \in \mathbb{ Z}$, so $  N ( p, q, j, \ell) \le \delta  $ if and only if $  p^2 k_1 + q^2 k_2 + j^2 + \ell^2 =0.$ Combined with conservation of the momentum, this gives \\
                  for the case $k = (-1,-1)$
                                    $$ p+q = j+ \ell \quad \text{and} \quad p^2 + q^2 = j^2 + \ell^2$$
                  for the case $k = (-2,0)$
                                 $$ 2p = j+ \ell \quad \text{and} \quad2p^2  = j^2 + \ell^2$$
                  for the case $k = (0,-2)$
                                 $$ 2q = j+ \ell \quad \text{and} \quad 2q^2 = j^2 + \ell^2$$
                  for the case $k = (-3,1)$
                                 $$ 3p = q+ j+ \ell \quad \text{and} \quad3p^2  = q^2+ j^2 + \ell^2$$
                  for the case $k = (1,-3)$
                                 $$ 3q = p+j+ \ell \quad \text{and} \quad 3q^2 = p^2+ j^2 + \ell^2.$$               
                  In all these cases, we get $ j,\ell \in \{ p,q\}$  which is excluded. 
     
       \textbf{The second case}  We see that $ \Omega$ and $\{\Lambda_j\}_{j \neq p,q,s,t}$ are all the same as the previous case except  $ \Lambda_{t^+}$ and $\Lambda_{t^-}$.
       We remind that 
      $$ \begin{cases}
      2p + s &= 2q +t\\
      2p^2 + s^2 &= 2q^2 +t^2.
      \end{cases}$$  
      Thank to Lemma 2.2 in \cite{quinticNLSeq}, $\{p,q,s,t\}$ is in form of $\{p,\,p+2n,\, p+3n, \,p-n\}.$ Without loss of generality, we can assume that\footnote{using the change of variables j=j-p} $p=0,$ so we have $q=-2t. $
      For $ \Omega \cdot k + \Lambda_{t^+}$ and  $ \Omega \cdot k + \Lambda_{t^-},$ by conservation the momentum, we just need to consider the case when $k$ satisfies  $p k_1 + q k_2 + t= 0$ i.e. $k_2=1/2,$ which is not an integer.  For $ \Omega \cdot k + \Lambda_{t^{\pm}} \pm \Lambda_j$, again by conservation of the momentum, we have 
                          $$ \begin{cases}
                                             p k_1 + q k_2 + t \pm j &= 0\\
                                             p^2 k_1 + q^2 k_2 + t^2 \pm j^2 &= 0
                                   \end{cases}$$
        i.e.
                          $$ \begin{cases}
                           j &= \mp ( 2k_2 -1)n\\
                          j^2 &= \mp(4k_2 +1)n^2.
                          \end{cases}$$
       This system has two solutions for $j$, either $j=0$($=p$) or $j=3m$($=s$), which are both excluded. 
                      
   \subsection{Three modes case.}
   It is too complicated to verify all the possibility, in this appendix we just do with an implicit example where $(p,q,m)=(-3,10,-6),$ which we are interesting in Theorem \ref{thm3}. In this situation, we have $\mathcal{C}, \, \mathcal{E}$ are all empty, $ \mathcal{A}=\{ -14, 2\}$ and  $\mathcal{B}=\{9,1\}$. 
   Recall that
              	$$ \Omega ( \rho) = 
              \begin{pmatrix}
              p^2 + 3 \nu^2 \left( \rho_1^2 + 3 \rho_2^2 + 3 \rho_3^2 +   6 \rho_1 \rho_2 + 6 \rho_1 \rho_3 +12 \rho_2 \rho_3\right)\\
              q^2 + 3 \nu^2 \left( \rho_2^2 + 3 \rho_1^2 + 3 \rho_3^2 +   6 \rho_1 \rho_2 + 6 \rho_2 \rho_3 +12 \rho_1 \rho_3 \right)\\
              m^2 + 3 \nu^2 \left( \rho_3^2 + 3 \rho_1^2 + 3 \rho_2^2 +   6 \rho_1 \rho_3 + 6 \rho_3 \rho_2 +12 \rho_2 \rho_1\right)
              \end{pmatrix}$$
      and
              $$  \Lambda_{j} = j^2 + 9\nu^2 ( \rho_1^2  + \rho_2^2 +\rho_3^2  + 4 \rho_1 \rho_2 + 4 \rho_2 \rho_3 + 4 \rho_3 \rho_1 ) \qquad j\neq -14,-6,-3,2,1,9,10.$$
        The hypothesis $A0$ and $A1$ are trivial. For hypothesis A2 $(1)$, let $k = ( k_1, k_2, k_3) \in \mathbb{Z }^3/ \{ 0\}$, $k'= (k_2+k_3, k_1+k_3,k_2+k_1)$ and $ z = z( k) = \frac{k'}{ | k'|},$ then we have
                   \begin{align*}  
                  	( \nabla_{\rho} \cdot z) ( \Omega (\rho) \cdot k) = &6\nu^2|k'|^{-1}[ \rho_1( 3k_2^2+3k_3^2+k_1k_2 +k_1k_3+ 6(k_1+k_2+k_3)^2 ) \\&+ \rho_2( 3k_1^2+3k_3^2+k_1k_2 +k_2k_3+  6(k_1+k_2+k_3)^2 )\\&+\rho_3( 3k_2^2+3k_1^2+k_3k_2 +k_1k_3+ 6(k_1+k_2+k_3)^2 ) ].  
                  \end{align*}
        This term is greater than $ \delta = \nu^2 $. Since $(\nabla_{\rho} \cdot z) (\Lambda_{j} - \Lambda_{\ell} ) =0 ,$ the estimate of small divisor $\Omega \cdot k + \Lambda_{j} - \Lambda_{\ell} $  is followed.\\
        For hypothesis A2 $(2),(3),$ choose $ z = z( k) = -\frac{k}{ | k|},$ then we have 
       \begin{align*}  
       	( \nabla_{\rho} \cdot z) ( \Omega (\rho) \cdot k) = &-6\nu^2|k|^{-1}[ \rho_1( k_1^2+3k_2^2+3k_3^2+6k_1k_2 +6k_1k_3+12k_2k_3 ) \\&+ \rho_2( k_2^2+3k_1^2+3k_3^2+6k_1k_2 +6k_2k_3+12k_1k_3 )\\&+\rho_3( k_3^2+3k_2^2+3k_1^2+6k_3k_2 +6k_1k_3+12k_2k_1 ) ]  
       \end{align*}
       and 
       $$ (\nabla_{\rho} \cdot z) \Lambda_{j}  = -18 \nu^2 | k|^{-1}[ \rho_1(k_1 +2k_2 +2k_3)+\rho_2(k_2 +2k_1 +2k_3)+\rho_3(k_3 +2k_2 +2k_1) ]
       .$$
         For $ \Omega \cdot k + \Lambda_{j},$ by conservation of the mass, we just need to estimate this divisor in the case $ k_1 + k_2 + k_3 =-1,$ then by computation we have 
         \begin{align*}
        	 	 | ( \nabla_{\rho} \cdot z) ( \Omega (\rho) \cdot k + \Lambda_j) &= 6 \nu^2 |k|^{-1} [ \rho_1 ( 2k_1^2 - 6k_2k_3 + 3k_1 +3 ) +  \rho_2 ( 2k_2^2 - 6k_1k_3 + 3k_2 +3 )\\&+ \rho_3 ( 2k_3^2 - 6k_2k_1 + 3k_3 +3 )]\\&\ge6 \nu^2 |k|^{-1} [ \rho_1 ( 2k_1^2 - \frac{3}{2}(k_1+1)^2 + 3k_1 +3 ) +  \rho_2 ( 2k_2^2 \\&- \frac{3}{2}(k_2+1)^2 + 3k_2 +3 )+ \rho_3 ( 2k_3^2 - \frac{3}{2}(k_3+1)^2 + 3k_3 +3 )]\\&= 3 \nu^2 |k|^{-1}[ \rho_1 ( k_1^2 + 3) + \rho_2 ( k_2^2 + 3) +\rho_3 ( k_3^2 + 3)  ] \\ &\ge \nu^2.
        	\end{align*}
        For   $\Omega \cdot k + \Lambda_{j} +\Lambda_{\ell}, $ again we have $k_1 + k_2 + k_3 =-2 $ by conservation of the mass, hence
          \begin{align*}
        	| ( \nabla_{\rho} \cdot z) ( \Omega (\rho) \cdot k + \Lambda_j) = &6 \nu^2 |k|^{-1}  [ \rho_1 ( 2k_1^2 - 6k_2k_3 + 6k_1 +12 ) +  \rho_2 ( 2k_2^2 - 6k_1k_3 + 6k_2 +12 )\\&+ \rho_3 ( 2k_3^2 - 6k_2k_1 + 6k_3 +12)]\\&\ge6 \nu^2 |k|^{-1} [ \rho_1 ( 2k_1^2 - \frac{3}{2}(k_1+1)^2 + 6k_1 +12 ) +  \rho_2 ( 2k_2^2 \\&- \frac{3}{2}(k_2+2)^2 + 6k_2 +12)+ \rho_3 ( 2k_3^2 - \frac{3}{2}(k_3+2)^2 + 6k_3 +12 )]\\&= 3 \nu^2 |k|^{-1}[ \rho_1 ( k_1^2 + 12) + \rho_2 ( k_2^2 + 12) +\rho_3 ( k_3^2 + 12)  ] \\ &\ge \nu^2.
        \end{align*}

        \textbf{The set $\mathcal{ B}$}  For $ \rho \in \mathcal{D}_2$: we have $$ | \Im \Lambda_{1^{\pm}}| > \nu^2= \delta$$ so that 
                                            $$ |\Omega \cdot k + \Lambda_{1^+} - \Lambda_{1^-} | \ge 2 \nu^2> \delta.$$
        For $\Omega \cdot k + \Lambda_{1^+} + \Lambda_{1^-},$ 
         by the conservation of the mass and the momentum, we just need to estimate this small divisor if
                 $$ \begin{cases}
                   k_1+ k_2 +k_3 +2&=0\\
                   -3k_1+ 10k_2 -6k_3 +2&=0\\
                   9k_1+ 100k_2 +36k_3 +2&=0\\
                   k_1,k_2,k_3 \in \mathbb{Z}
                  \end{cases}	$$
        This equation system has no solution\footnote{with the implicit form of $\{p,q,m,s,t\}$ in appendix B, we can solve for general $p,q,m$}. 
        
        \textbf{The set $\mathcal{ A}$} For $ \Omega \cdot k + \Lambda_{2^{\pm}}$ and $ \Omega \cdot k + \Lambda_{2^{\pm}} + \Lambda_j $ again by the conservation of the mass and the momentum, we have   
        \begin{align*}
        (*)  \begin{cases}
        k_1+ k_2 +k_3 +1&=0\\
        -3k_1+ 10k_2 -6k_3 +2&=0\\
        9k_1+ 100k_2 +36k_3 +4&=0
        \end{cases}	
        \qquad (**) 
          \begin{cases}
        k_1+ k_2 +k_3 +2&=0  \\
        -3k_1+ 10k_2 -6k_3 +2 +j&=0 \\
        9k_1+ 100k_2 +36k_3 +4+ j^2&=0.
        \end{cases}
        \end{align*}
      It is easy to see that $(*)$ has no solution in $\mathbb{Z}^3.$ For $(**)$ we have $ j \equiv -k_2 -2 \pmod{3} $ and $ j^2 \equiv -k_2 -4 \pmod{9} $. If $j \equiv \pm 1 \pmod{3}$ then we have $ k_2 \equiv 0,2 \pmod{4} $ and $k_2=4 \pmod{9},$ which can not both happen. If $j \equiv 0 \pmod{3}$ then we have $ k_2 \equiv 1 \pmod{4} $ and $k_2=5 \pmod{9},$ which again can not happen.  For $ \Omega \cdot k + \Lambda_{2^{\pm}} - \Lambda_j,$ because of changes of variables, we have 
  \begin{align*}
           \Lambda_{2^+}&= \Lambda_2 -g(\rho_1,\rho_2,\rho_3)\\
            \Lambda_{2^-}&= \Lambda_2 -g(\rho_1,\rho_2,\rho_3) +12(\rho_3^2-\rho_2^2+3\rho_1 \rho_3=3\rho_2\rho_1)
   \end{align*}           
 with $g(x,y,z)= \mu^2 \sqrt{81 y^2z^2 +(-18xy+18xz-6y^2+6z^2)^2}- \mu^2(-18xy+18xz-6y^2+6z^2).$
      By the conservation of the mass we just need to consider the case $ k_1+k_2+k_3=0,$ then 
      \begin{align*}
      (\nabla_{\rho} \cdot z) ( \Omega \cdot k + \Lambda_{2^{\pm}}- \Lambda_{j})=&12 \mu^2|k|^{-1} [\rho_1( k_2^2 + k_3^2 -k_2k_3 -2k_2 +k_3 )  \\&+ \rho_2( k_1^2 + k_3^2 -k_1k_3 +2k_1 -k_3 ) \\&+\rho_3( k_2^2 + k_1^2 -k_2k_1 +3k_1 -3k_2 )  ] \pm (\nabla_{\rho} \cdot z)g \\&\approx 12 |k| \mu^2 \rho \pm | (\nabla_{\rho} \cdot z)g |.
     \end{align*} 
     By the conservation of the momentum we have 
                  \begin{align*}   
                  \begin{cases}
                  -3k_1+ 10k_2 -6k_3 +2 -j&=0 \\
                  9k_1+ 100k_2 +36k_3 +4- j^2&=0.
                  \end{cases}                     
                  \end{align*} 
     The solution of this equation system that closest to the origin is $k=(-975,195,780)$ and with such a big $k,$ $  (\nabla_{\rho} \cdot z) ( \Omega \cdot k + \Lambda_{2^{\pm}}- \Lambda_{j})$ is far greater than $\delta.$



       
 \section{Appendix B}
      In this appendix, we try to solve the set $\mathcal{ B}$ in general
                     $$ 
                     \begin{cases}
                      2p+ q&=m+s+t\\
                      2p^2 + q^2 &= m^2+ s^2+ t^2.
                     \end{cases}     
                     $$
     Let $q'=q-p,\, m'=m-p,\, s'=s-p,\, t'=t-p,$ it becomes
                    $$ 
                   \begin{cases}
                    q'&=m'+s'+t'\\
                   q'^2 &= m'^2+ s'^2+ t'^2.
                   \end{cases} 
                   $$     
     This give us $ m's'+ t's'+ t'm'=0,$ hence $ s'=- \frac{m't'}{m'+t'}.$ Assume more that $s',t',m'$ have no common divisor except $\pm 1.$ Let $k$ is a prime common divisor of $t'$ and $m',$ i.e. $t'=t" k,\, m'=m"k,$ then $ s'=- \frac{km"t"}{ m"+t"}.$ Since $ k\nmid s,$ we have $ k \mid t" +m",$ i.e. $t"= kh-m",$ hence $ s'= -\frac{m"(kh-m")}{h}=-km"+ \frac{m"^2}{h} \in \mathbb{Z }.$ Let $h=(-1)^{sgn(h)} \Pi p_i^{k_i},$ $x= \Pi p_i^{[\frac{k_i}{2}]}$ and $y= (-1)^{sgn(h)}\Pi p_i^{ k_i -2 [\frac{k_i}{2}]},$ with $p_i$ is prime divisor of $h.$ Then, $h= x^2y$ and we need $xy \mid m",$ i.e. $m"= ryx.$ By this, $ s' = -kxyr + r^2y,\, m'= kryx,\, t'= k^2x^2y - ryx.$ Since $s',t',m'$ have no common divisor except $\pm 1,$ we have $y= \pm 1.$ Assume that $y=1,$ and $kx=n,$ then $s'= r^2 -nr,\, m'=nr,\, t'=n^2-nr$ and $q'=n^2-nr+r^2.$ In general, we have $\{p,q,m,s,t\}= \{p, p+ k(n^2-nr+r^2),p+ knr,p+ k(r^2-nr),p+ k(n^2-nr)\}.$

\end{document}